\documentclass{amsart}
\newcommand\str{{\text{\raisebox{0.5ex}{\normalfont\bf\Large.}}}}
\usepackage
{xypic}
\numberwithin{equation}{section}
\newtheorem{Lem}{Lemma}[section]
\newtheorem{Prop}[Lem]{Proposition}
\newtheorem{Cor}[Lem]{Corollary}
\newtheorem{Thm}[Lem]{Theorem}

\theoremstyle{definition}

\theoremstyle{remark}
\newtheorem{Rem}[Lem]{Remark}

\newcommand{\mapsfrom}{\mathrel{\mbox{$\leftarrow\joinrel\mapstochar$\,}}}
\renewcommand\o{\otimes}
\DeclareMathOperator\Hom{\operatorname{Hom}}
\DeclareMathOperator\End{\operatorname{End}}

\DeclareMathOperator\ev{\operatorname{ev}}
\DeclareMathOperator\db{\operatorname{db}}

\newcommand\bop{{\operatorname{bop}}}

\newcommand\drdot[1]{{^\str_\str #1^\str}}
\newcommand\HMod[4]{{^{#1}_{#3}\mathcal M^{#2}_{#4}}}
\newcommand\LMod[1]{{_{#1}\mathcal M}}
\newcommand\LComod[1]{{^{#1}\mathcal M}}
\newcommand\RComod[1]{\mathcal M^{#1}}
\newcommand\RComodf[1]{\mathcal M^{#1}_f}
\newcommand\RMod[1]{\mathcal M_{#1}}
\newcommand\LModf[1]{{_{#1}\mathcal M_f}}
\newcommand\BiMod[1]{{_{#1}\mathcal M_{#1}}}
\newcommand\BiComod[1]{{^{#1}\mathcal M^{#1}}}

\newcommand\kmod{\mathcal M_k}

\newcommand{\leer}{\operatorname{--}}
\newcommand{\ou}[1]{\mathrel{\mathop{\otimes}_{#1}}}
\newcommand{\co}[1]{\mathrel{\mathop{\Box}_{#1}}}

\newcommand\si[1]{(#1)}
\newcommand\sw[1]{{}_{(#1)}}
\newcommand\swm[1]{{}_{(-#1)}}
\newcommand\so[1]{^{(#1)}}
\newcommand\som[1]{^{(-#1)}}

\newcommand\ol{\overline}

\newcommand\inv{^{-1}}

\renewcommand\epsilon\varepsilon

\newcommand\xihut{\hat\xi}
\newcommand\pEv{E'}
\newcommand\Ev{E}
\newcommand\Db{D}
\newcommand\tgamma{\widetilde\gamma}
\newcommand\ttheta{\widetilde\vartheta}
\newcommand\tttheta{\vartheta^t}
\newcommand\imp{\ensuremath{\implies}}
\newcommand\C{\mathcal C}
\newcommand\dual[1]{#1^\vee}
\newcommand\bidual[1]{#1^{\vee\vee}}
\newcommand\catr{\underline{\operatorname{tr}}}
\newcommand\tr{\operatorname{tr}}
\newcommand\piv{j}
\def\Hfix#1.{#1^H}
\newcommand\rdual[1]{{{\,^\vee\! #1}}}
\newcommand\li{^{(1)}}
\newcommand\re{^{(2)}}
\newcommand\FSI{\nu_V}
\makeatletter
\def\namelabel#1#2{\@bsphack
  \protected@write\@auxout{}%
         {\string\newlabel{#1.nme}{{#2}{#2}}}%
  \@esphack}
\def\nmlabel#1#2{\label{#2}\namelabel{#2}{#1}}
\newcommand\nmref[1]{\ref{#1.nme}\ \ref{#1}}
\makeatother

\begin{document}
\title{On the Frobenius-Schur indicators for quasi-Hopf algebras}
\author{Peter Schauenburg}
\address{Department of Mathematics, USC, Los Angeles, CA 90089-1113, USA, and Mathematisches Institut der Universit\"at M\"unchen,
Theresienstr.~39, 80333~M\"unchen, Germany}
\email{schauen@mathematik.uni-muenchen.de} 
\thanks{The author is supported by a Heisenberg-Stipendium from the DFG (Deutsche Forschungsgemeinschaft)}
\begin{abstract}
Mason and Ng have given a generalization to semisimple quasi-Hopf
algebras of Linchenko and Montgomery's generalization to
semisimple Hopf algebras of the classical Frobenius-Schur theorem
for group representations. We give a simplified proof, in
particular a somewhat conceptual derivation of the appropriate
form of the Frobenius-Schur indicator that indicates if and in
which of two possible fashions a given simple module is self-dual.
\end{abstract}

\maketitle
\section{Introduction}
Let (throughout this paper) $k$ be an algebraically closed field
of characteristic zero.

Let $G$ be a finite group, and $V$ a finite-dimensional $k$-linear
representation of $G$. According to the classical Frobenius-Schur
Theorem, $V$ falls into one of the following three mutually
exclusive classes: Either it is not selfdual, or there is a
nondegenerate $G$-invariant bilinear form on $V$ that is
symmetric, respectively skew-symmetric. Which of the three cases
is at hand can be decided by computing the Frobenius-Schur
indicator $\nu_V=\chi_V(|G|\inv\sum_{g\in G}g^2)$. This can only
take the values $0,1,-1$, which occur when $V$ is not selfdual,
admits a nondegenerate symmetric invariant bilinear form, or
admits a nondegenerate skew-symmetric invariant bilinear form,
respectively.

Linchenko and Montgomery \cite{LinMon:FSTHA} have generalized the
Frobenius-Schur Theorem to irreducible representations $V$ of a
semisimple Hopf algebra $H$, which are found to fall into the same
three classes, again detected by the Frobenius-Schur indicator,
now given by $\nu_V=\chi_V(\Lambda\sw 1\Lambda\sw 2)$ for a
normalized integral $\Lambda\in H$.

Finally, Mason and Ng \cite{MasNg:CIFSISQHA} have generalized the
results of Linchenko and Montgomery to the case of a semisimple
quasi-Hopf algebra $H$. As pointed out by Mason and Ng (and also
by Pavel Etingof in an email to Susan Montgomery), a quite general
argument essentially due to Fuchs, Ganchev, Szl\'achanyi, and
Vescerny\'es \cite{FucGanSzlVes:FSIRMCC} shows that a trichotomy
generalizing that in the Frobenius-Schur and Linchenko-Montgomery
Theorems can be obtained for simple objects of any pivotal fusion
category --- this covers the representation category of a
semisimple quasi-Hopf algebra by  results of Etingof, Nikshych,
and Ostrik \cite{EtiNikOst:FC}. It is perhaps more surprising that
a formula for a Frobenius-Schur indicator detecting the three
cases can again be found, by adorning the formula from
\cite{LinMon:FSTHA} with the structure elements of a quasi-Hopf
algebra in a rather non-obvious fashion.

The present paper contains no new results in this direction. We
do, however, present a new proof for the fact that the
Frobenius-Schur indicator of Mason and Ng does in fact indicate
what it should. We find this worthwhile because our proof does not
involve a lot of calculations with the structure elements and
axioms of quasi-Hopf algebras, and is quite short. We should
acknowledge that on our short and direct route to the indicator
formula we do not encounter (and prove) many of the facts
contained in the work of Mason and Ng.

Key to our proof is writing the Frobenius-Schur indicator of a
module $V$ as the trace of a certain endomorphism of the invariant
subspace of $V\o V$. In the case of ordinary Hopf algebras, where
the endomorphism is given by flipping tensor factors, this
argument appears in work of Kashina, Sommerh\"auser, and Zhu
\cite{KasSomZhu:HFSI}.

\section{Definitions and Preliminaries}
A monoidal category $\C$ is a category equipped with a bifunctor
$\o\colon\C\times\C\to\C$, a coherent natural (in $U,V,W\in\C$
isomorphism $\Phi\colon (U\o V)\o W\to U\o(V\o W)$, and a neutral
object $I$, which, for simplicity, we assume strict in the sense
that $I\o V=V=V\o I$, naturally in $V\in\C$. Coherence of $\Phi$
(and the neutral object) means that all reasonable diagrams
involving only instances of $\Phi$ and introducing and cancelling
neutral objects commute. A consequence is that $\C$ is equivalent
as a monoidal category to a strict monoidal category, that is, one
in which $\Phi$ is an identity. We refer to \cite{Kas:QG} for
background on monoidal categories.

Except for a few facts on general monoidal categories, we will
mostly deal with vector spaces over a fixed algebraically closed
field $k$ of characteristic zero, and $\o$ will mean tensor
product over $k$. We let $\tau\colon V\o W\to W\o V$ denote the
flip of tensor factors. We will recklessly use the formal notation
$v\o w\in V\o W$ for general elements of $V\o W$, even if we know
that the element in question is not a simple tensor. In
particular, notations like $q=q\so 1\o q\so 2\in V\o W$ will often
be used to refer to the components of a tensor.

A quasibialgebra $H$ is an associative algebra (whose
multiplication map we denote by $\nabla\colon H\o H\to H$)
with an
algebra map $\Delta\colon H\rightarrow H\o H$ called
comultiplication, an algebra map $\epsilon\colon H\rightarrow k$
that is a counit for $\Delta$, and an invertible element $\phi\in
H\o H\o H$ such that the category $\LMod H$ is a monoidal category
as follows: The tensor product of $V,W\in\LMod H$ is taken over
$k$ and endowed with the diagonal module structure $h(v\o w)=h\sw
1v\o h\sw 2w$, where $\Delta(h)=h\sw 1\o h\sw 2$, the coherent
associator isomorphism $\Phi\colon (U\o V)\o W\to U\o (V\o W)$ for
$U,V,W\in\LMod H$ is left multiplication by $\phi$. We write
$\phi=\phi\so 1\o\phi\so 2\o\phi\so 3$ and $\phi\inv=\phi\som
1\o\phi\som 2\o\phi\som 3$ when in need of referring to the tensor
components of $\phi$. The neutral object in $\LMod H$ is the
trivial module $k$ with action via $\epsilon$, and the usual
canonical identifications $k\o V=V=V\o k$ for $V\in\LMod H$.

We made a point of not actually recalling the axioms on the data
$\Delta$ and $\phi$ that ensure that $\LMod H$ is monoidal --- we
will never have to use most of those axioms. We will need to know
that applying $\epsilon$ to any of the three tensor factors of
$\phi$ yields $1\in H\o H$. For details and background on
quasibialgebras and quasi-Hopf algebras, which were introduced by
Drinfeld \cite{Dri:QHA}, we again refer to \cite{Kas:QG}.

A quasiantipode $(S,\alpha,\beta)$ for a quasibialgebra $H$
consists of an anti-algebra endomorphism $S$ of $H$, and elements
$\alpha,\beta\in H$, such that
\begin{align*}
  S(h\sw 1)\alpha h\sw 2&=\epsilon(h)\alpha, & h\sw 1\beta S(h\sw 2)&=\epsilon(h)\beta,\\
  \phi\so 1\beta S(\phi\so 2)\alpha\phi\so 3&=1,&S(\phi\som 1)\alpha\phi\som 2\beta\phi\som 3&=1
\end{align*}
hold in $H$, for $h\in H$. A quasi-Hopf algebra is a
quasibialgebra with a quasi-antipode. If $H$ is finite-dimensional
then $S$ is always a bijection by a result of Bulacu and Caenepeel
\cite{BulCae:IDQHAA}. If $H$ has a quasiantipode, then the
category of finite dimensional $H$-modules is rigid: Every
finite-dimensional module $V\in\LMod H$ has a dual object $\dual
V\in\LMod H$, whose underlying vector space is the vector space
dual $V^*$, endowed with the module structure in which $h\in H$
acts as the dual of the map given by the action of $S(h)$ on $V$.
The evaluation and coevaluation (or ``dual basis'') morphisms are
given by
  \begin{gather*}
    \ev\colon V^*\o V\ni\varphi\o v\mapsto \varphi(\alpha v)\in k\\
    \db\colon k\ni 1\mapsto \beta v_i\o v^i.
  \end{gather*}
where $v_i\o v^i\in V\o V^*$ is the canonical element.

If $H$ has bijective antipode, then the functor $\dual{(\leer)}$
is an equivalence. We denote the inverse functor by
$\rdual{(\leer)}$ and call $\rdual V$ the right dual of $V$. Just
like the left dual, the right dual has $V^*$ as its underlying
vector space.

As for any object $V$ of a monoidal category having a dual object
$\dual V$, we obtain an adjunction
$$A_0\colon \Hom_H(U,V\o W)\cong\Hom_H(\dual V\o U,W)$$
by
\begin{gather*}
  A_0(f)=\left(\dual V\o U\xrightarrow{\dual V\o f}\dual V\o (V\o
  W)\xrightarrow{\Phi\inv}(\dual V\o V)\o W\xrightarrow{\ev\o
  W}W\right)\\
  A_0\inv(g)=\left(U\xrightarrow{\db\o U}(V\o\dual V)\o
  U\xrightarrow\Phi V\o(\dual V\o U)\xrightarrow{V\o g}V\o
  W\right)
\end{gather*}
We recall the following elements of $H\o H$ introduced by Hausser
and Nill \cite{HauNil:DCPDQQG,HauNil:DQQG}:
\begin{align*}
  q_R&=\phi\so 1\o S\inv(\alpha\phi\so 3)\phi\so 2 &p_R&=\phi\som
  1\o\phi\som 2\beta S(\phi\som 3)\\
  q_L&=S(\phi\som 1)\alpha\phi\som 2\o\phi\som 3&
  p_L&=\phi\so 2S\inv(\phi\so 1\beta)\o\phi\so 3
\end{align*}
Writing $f(u)=f\so 1(u)\o f\so 2(u)\in V\o W$ for $f\colon U\to
V\o W$, we find
$$A_0(f)(\varphi\o u)=(\phi\som 1\varphi)(\alpha\phi\som 2 f\so
1(u))\phi\som 3f\so 2(u)=\varphi(q_L\so 1f\som 1(u))q_L\so 2f\som
2(u)
$$ and
$$A_0\inv(g)(u)=\phi\so 1\beta v_ig(\phi\so 2 v^i\o \phi\so 3u)
=v_ig(p_L\so 1v^i\o p_L\so 2u).$$ Specializing $U=k$ and using the
identification $\Hom_H(k,M)=M^H=\{m\in M|hm=\epsilon(h)m\forall
h\in H\}$, we find
\begin{equation}\label{A}
A\colon (V\o W)^H\overset\sim\rightarrow \Hom_H(\dual V,W);\quad
A(v\o w)(\varphi)=\varphi(q_L\so 1v) q_L\so 2w
\end{equation}
with $A\inv(g)=\beta v_ig(v^i)$. We note $v_i\o A(v\o
w)(v^i)=q_L(v\o w)$.

Finally, recall \cite{HauNil:ITQHA} that a quasi-Hopf algebra $H$
is semisimple if and only if it contains a normalized integral
$\Lambda$, that is, an element $\Lambda\in H$ satisfying
$h\Lambda=\epsilon(h)\Lambda$ for all $h\in H$ and
$\epsilon(\Lambda)=1$. For any $H$-module $M$, the $H$-linear
projection $\pi\colon M\to M^H$ onto the isotypical component of
the trivial module, or the invariant subspace $M^H=\{m\in
M|hm=\epsilon(h)m\forall h\in H\}$, is multiplication by the
normalized integral $\Lambda$.

In the sequel, we will refer to all the structure elements
mentioned for a (semisimple) quasi-Hopf algebra in this section
without further notice.

\section{The pivotal structure}
A pivotal structure for a monoidal category $\C$ with a duality
functor $\dual{(\leer)}$ is by definition a monoidal natural
isomorphism $\piv\colon V\to \bidual V$. We note that (as a
special case of a fact pertaining to all monoidal transformations)
we have $\dual{(\piv_V)}=\piv_{\dual V}\inv$ for all $V\in\C$.

Recall \cite{EtiNikOst:FC} that for a morphism $f\colon V\to
\bidual V$ in a rigid monoidal category one can define the
(categorical) trace of $f$ by
$$\catr_V(f)=\left(I\xrightarrow{\db}V\o \dual V\xrightarrow{f\o \dual V}\bidual
  V\o \dual V\xrightarrow{\ev}I\right)\in\End_\C(V)$$

We will depend heavily upon the following result of Etingof,
Nikshych, and Ostrik: By \cite[Prop.8.24 and
Prop.8.23]{EtiNikOst:FC}, the monoidal category $\LMod H$ for a
semisimple quasi-Hopf algebra is pivotal. More precisely, for any
simple $H$-module $V$ there is a unique isomorphism $\piv_V\colon
V\to \bidual V$ with $\catr_V(\piv_V)=\dim V$. The $\piv_V$ are
the components at the simple modules of a unique natural
transformation $\piv$, which is a pivotal structure on $\LMod H$,
and satisfies $\catr_V(\piv_V)=\dim V$ for all finite-dimensional
$H$-modules $V$.

Since $\piv$ is natural, there exists a unique invertible $g\in H$
with $\piv(v)=\iota(gv)$ for all $v\in V$, where $\iota\colon V\to
V^{**}$ is the canonical vector space isomorphism.

Since the invariant subspace $\Hfix (V\o W).$ is isomorphic to
$\Hom_H(k,V\o W)$, we see that $\dual V\o V$ contains a unique up
to scaling $H$-invariant element for every simple $V$. The formula
for the categorical trace of $\piv_{\rdual V}$ shows that the
unique invariant element in $\dual V\o V$ mapped to $\dim V$ under
the evaluation map $\dual V\o V\to k$ is $(\piv_{\rdual V}\o
V)\db_{\rdual V}(1)=\piv(\beta v^i)\o v_i=\beta\piv(v^i)\o
v_i=\beta v^i\o\piv\inv(\iota(v_i))=\beta v^i\o g\inv v_i$.

On the other hand, it is easy to exhibit elements in $\dual V\o V$
mapped to $\dim V$ under the evaluation map: For $t=t\li\o t\re\in
H\o H$, the element $t\li v^i\o t\re v_i$ is mapped to $(t\li
v^i)(\alpha t\re v_i)=\chi_V(S(t\li)\alpha t\re)$. For this to be
equal to $\dim V=\chi_V(1)$ it suffices to have $S(t\li)\alpha
t\re=1$, which is the case both for $t=p_L$ and $t=p_R$. Since the
evaluation morphism is $H$-linear, we have $\ev(\Lambda
\xi)=\ev(\xi)$ for every $\xi\in\dual V\o V$. If $\ev(\xi)=\dim
V$, then it follows that $\Lambda\xi=\beta v^i\o g\inv v_i$. We
conclude
\begin{Lem}
  Let $H$ be a semisimple quasi-Hopf algebra and $t\in\{p_L,p_R\}\subset H\o H$.

  For every finite-dimensional $H$-module $V$ we have
  $$\beta v^i\o g\inv v_i=\Lambda\sw 1t\li v^i\o\Lambda\sw 2t\re v_i\in \dual V\o V$$
  and $g\inv S(\beta)=S(\Lambda\sw 1t\li)\Lambda\sw 2t\re\in H$.
\end{Lem}
\begin{proof}
  We proved the first formula for the case where $V$ is simple.
  This entails that the second formula holds in $\End_k(V)$ for
  every simple module, hence in $H$, which in turn implies the
  first formula for any module.
\end{proof}
\begin{Rem}
  The second formula gives an explicit formula for $g$ in the case
  that $\beta$ is a unit in $H$, and comes close otherwise.
  However, the first, somewhat less explicit formula proves
  particularly useful in our computations below.
\end{Rem}

\section{Frobenius-Schur indicators}
Let $V$ and $W$ be two objects of a pivotal rigid monoidal
category $\C$. Then we have a bijection
$$T=T_{VW}\colon\C(\dual V,W)\to\C(\dual W,V),$$
given by $T(f)=(\dual W\xrightarrow{\dual f}\bidual
V\xrightarrow{\piv_V\inv}V)$. We note that
$$  T_{WV}T_{VW}(f)
  =j_W\inv\dual{(j_V\inv\dual f)}=j_W\inv\bidual
  f\dual{(j_V\inv)}=fj_{\dual V}\inv \dual{(j_V\inv)}=f
$$
In particular, $T_{VV}$ is an order two automorphism of $\C(\dual
V,V)$. In the case that $V$ is a self-dual simple module in
$\C=\LMod H$, $T_{VV}$ is multiplication by a scalar thanks to
Schur's lemma, and the scalar can only be $\pm 1$. Thus, there are
exactly three possibilities for a simple $H$-module $V$: Either
$V$ is not self-dual, in which case we put $\mu_V=0$, or it is,
and every isomorphism $f\colon \dual V\to V$ satisfies
$\piv_Vf=\mu_V \dual f$, with the sign $\mu_V\in\{\pm 1\}$
independent of the choice of $f$.

Linchenko and Montgomery \cite{LinMon:FSTHA}, generalizing a
classical result on group representations by Frobenius and Schur,
have proved the above trichotomy for simple modules of a
semisimple Hopf algebra (note $\piv$ is trivial in this case by
the Larson-Radford Theorem), and shown that which of the cases
applies to a particular module can be decided using the
Frobenius-Schur indicator, given by the formula
$\FSI(V)=\chi_V(\Lambda\sw 1\Lambda\sw 2)$. The indicator can only
take the values $0,\pm 1$, depending on which of the three
possible behaviors $V$ displays.

Mason and Ng \cite{MasNg:CIFSISQHA} have in turn generalized the
results of Linchenko and Montgomery to the quasi-Hopf algebra
case. We have already seen that the same trichotomy is preserved
here. It can be reformulated as in \cite{MasNg:CIFSISQHA} to be a
property of bilinear forms on $V$: An isomorphism $f\colon V\to
V^*$ is the same as a nondegenerate bilinear form $\sigma\colon
V\o V\to k$, by the correspondence $\sigma(v,w)=f(v)(w)$. One
checks that $f$ is an $H$-linear map $f\colon V\to \dual V$ if and
only if $\sigma(hv,w)=\sigma(v,S(h)w)$, that is, $\sigma$ is a
bilinear form with adjoint $S$ in the sense of
\cite{MasNg:CIFSISQHA}. If $f$ corresponds to $\sigma$, then
$\dual f\piv_V$ corresponds to $\sigma'$ given by
$\sigma'(v,w)=\sigma(w,gv)$. Thus, the trichotomy can be restated
as follows: If there exists a nondegenerate bilinear form $\sigma$
on $V$ with adjoint $S$, then there is $\mu_V\in\{\pm 1\}$ so that
$\sigma(v,w)=\mu_V\sigma(w,gv)$ for any such form. If there is no
such form, we put $\mu_V=0$ again.

In the rest of this section we will deal with the question how one
can come up with a modification of the Frobenius-Schur indicator
of Linchenko and Montgomery suitable for the quasi-Hopf algebra
case. This means that, rather than only label the three possible
relations of $V$ to its dual by a number $\mu_V$, we are looking
for a formula computing a number $\nu_V\in k$ from the structure
elements of a semisimple quasi-Hopf algebra and the character
$\chi_V$ of a simple module $V$, such that $\nu_V=\mu_V$ (and in
particular $\nu_V\in\{0,\pm 1\}$).

We begin by considering two $H$-modules $V,W$ and computing the
composition
$$E_{VW}:=\left((V\o W)^H\xrightarrow A \Hom_H(\dual
V,W)\xrightarrow{T_{V,W}}\Hom_H(\dual W,V)\xrightarrow{A\inv}(W\o
V)^H\right),$$ where $A$ was introduced in \eqref{A}. We have
\begin{align*}
  A\inv TA(v\o w)&=\beta w_i\o TA(v\o w)(w^i)\\
  &=\beta w_i\o \piv_V\inv \dual{A(v\o w)}(w^i)\\
  &=\beta A(v\o w)(v^i)\o \piv_V\inv(\iota(v_i))\\
  &=A(v\o w)(\beta v^i)\o g\inv(v_i)\\
  &=A(v\o w)(\Lambda\sw 1t\li v^i)\o \Lambda\sw 2t\re v_i\\
  &=\Lambda\sw 1t\li A(v\o w)(v^i)\o\Lambda\sw 2t\re v_i\\
  &=\Lambda\sw 1t\li q_L\so 2 w\o \Lambda\sw 2t\re q_L\so 1v
\end{align*}
for $v\o w\in(V\o W)^H$, so that
$$E_{VW}=\left((V\o W)^H\hookrightarrow V\o
W\xrightarrow{q_L\cdot}V\o W\xrightarrow\tau W\o
V\xrightarrow{t\cdot}W\o V\xrightarrow\pi(W\o V)^H\right).$$ Now
consider the case where $V=W$ is simple. Then $(V\o V)^H$ is zero
unless $V$ is selfdual, and then $(V\o V)^H$ is one-dimensional.
In either case $T_{VV}$ is multiplication by $\mu_V\in\{0,\pm
1\}$. Consequently $E_{VV}$ is multiplication by the same scalar
$\mu_V$, and this is also the trace of $E_{VV}$. But $E_{VV}$ has
the same trace as
\begin{align*}
  E'&:=\left(V\o V\xrightarrow\pi(V\o
V)^H\hookrightarrow V\o V\xrightarrow{q_L\cdot}V\o
V\xrightarrow\tau V\o V\xrightarrow{t\cdot}V\o V\right)\\
&=\left(V\o V\xrightarrow{q_L\Delta(\Lambda)\cdot}V\o
V\xrightarrow\tau V\o V\xrightarrow{t\cdot}V\o V\right).
\end{align*}
Now using the general observation
\begin{equation}\label{twisttrick}\tr((f_1\o f_2)\circ\tau\circ(f_3\o
f_4))=\tr(f_3f_1f_4f_2) \end{equation}
 for endomorphisms
$f_1,\dots,f_4$ of $V$ we find that
\begin{align*}
  \tr(E_{VV})&=\tr(E')\\
    &=\chi_V(q_L\so 1\Lambda\sw 1t\so 1q_L\so 2\Lambda\sw 2t\so
    2)\\
    &=\chi_V(\nabla(q_L\Delta(\Lambda)t)).
\end{align*}
We have proved the following result of Mason and Ng:
\begin{Thm}
  Let $H$ be a semisimple quasi-Hopf algebra, and
  $t\in\{p_L,p_R\}$. Put
  $\nu_H:=\nabla(q_L\Delta(\Lambda)t)$. Let $V$ be a simple
  $H$-module, and define the Frobenius-Schur indicator
  $\nu_V:=\chi_V(\nu_H)$. Then $\nu_V\in\{0,\pm 1\}$, and the following
  are equivalent for $\mu\in\{\pm 1\}$:
  \begin{enumerate}
    \item $\nu_V=\mu$
    \item There is an $H$-linear isomorphism $f\colon V\to\dual V$
    with $\dual f\piv_V=\mu f$.
    \item There is a nondegenerate bilinear form $\sigma$ on $V$
    with adjoint $S$ and $\sigma(v,w)=\mu\sigma(w,gv)$ for all
    $v,w\in V$.
  \end{enumerate}
\end{Thm}

In closing we should stress that the key trick of writing the
indicator as the trace of an endomorphism of $(V\o V)^H$, using
\eqref{twisttrick}, is due to Kashina, Sommerh\"auser, and Zhu
\cite{KasSomZhu:HFSI}; they treat the case of ordinary Hopf
algebras, where $E_{VW}$ is just the flip of tensor factors.

\bibliographystyle{acm}
\bibliography{promo}
\end{document}